\documentclass[8pt]{article}
\usepackage{amsfonts}
\usepackage{latexsym,amsmath}
\usepackage{amssymb,array}
\usepackage{graphics}
\makeatletter \oddsidemargin 0in \evensidemargin 0in \textwidth
16cm \RequirePackage[dvips]{graphicx} \textheight 20cm
\newtheorem{d1}{Definition}[section]

\newtheorem{ex}{Example}[section]

\begin{document}
\title{A Generalized Xgamma Generator Family of Distributions}
\author{Sudhansu S. Maiti$^1$\footnote{Corresponding author. e-mail:}, Sukanta Pramanik$^2$\\
$^1$Department of Statistics, Visva-Bharati University, Santiniketan-731 235, India\\
$^2$Department of Statistics, Siliguri College, North Bengal University, Siliguri, India}
\date{}
\maketitle

\begin{abstract}
In this paper, a new class of distributions, called Odds xgamma-G (OXG-G) family of distributions is proposed for modeling lifetime data. A comprehensive account of the mathematical properties of the new class including estimation issue is presented. Two data sets have been analyzed to illustrate its applicability.
\end{abstract}
\section{Introduction}
There are several ways of adding one or more parameters to a distribution function. Such an addition of parameters makes the resulting distribution richer and more flexible for modeling data. Proportional hazard model (PHM), Proportional reversed hazard model (PRHM), Proportional odds model (POM), Power transformed model (PTM) are few such models originated from this idea to add a shape parameter. In these models, a few pioneering works are by Box and Cox (1964),Cox (1972), Mudholkar and Srivastava (1993), Shaked and Shantikumar (1994), Marshall and Olkin (1997), Gupta and Kundu (1999), Gupta and Gupta (2007) among others.\\
Many distributions have been developed in recent years that involves the logit of the beta distribution. Under this generalized class of beta distribution scheme, the cumulative distribution function (cdf) for this class of distributions for the random variable $X$ is generated by applying the inverse of the cdf of $X$ to a beta distributed random variable to obtain,
\begin{eqnarray*}
	F(x)= \frac{1}{B(\alpha,~\beta)}\int_{0}^{G(x)}t^{\alpha-1}((1-t)^{\beta-1}dt;~\alpha,~\beta>0,
\end{eqnarray*}
where $G(x)$ is the cdf of any other distribution. This class has not only generalized the beta distribution but also added parameter(s) to it. Among this class of distributions are, the beta-Normal [Eugene et el. (2002)]; beta-Gumbel [Nadarajah and Kotz (2004)]; beta-Exponential [Nadarajah and Kotz (2006)]; beta-Weibull [Famoye et al. (2005)]; beta-Rayleigh [Akinsete and Lowe (2009)]; beta-Laplace [Kozubowski and Nadarajah (2008)]; and beta-Pareto [Akinsete et al. (2008)], among a few others. Many useful statistical properties arising from these distributions and their applications to real life data have been discussed in the literature.\\
In the generalized class of beta distribution, since the beta random variable lies between $0$ and $1$, and the distribution function also lies between $0$ and $1$, to find out cdf of generalized distribution, the upper limit is replaced by cdf of the generalized distribution.\\
Alzaatreh et al. (2013) has proposed a new generalized family of distributions, called T-X family, and the cumulative distribution function (cdf) is defined as
\begin{eqnarray}\label{eq1.1}
	F(x; \theta)= \int_{a}^{W[G(x)]}r(t)dt,
\end{eqnarray}
where, the random variable $T\in [a,b]$, for $-\infty<a,b<\infty$ and $W[G(x)]$ be a function of the cdf $F(x)$ so that $W[G(x)]$ satisfies the following conditions:
\begin{enumerate}
	\item[(i)] $W[G(x)]\in [a,b]$,
	\item[(ii)] $W[G(x)]$ is differentiable and monotonically non-decreasing,
	\item[(iii)] $W[G(x)]\rightarrow a$ as $x\rightarrow -\infty$ and $W[G(x)]\rightarrow b$ as $x\rightarrow \infty$.	
\end{enumerate}
In this paper, we propose a new wider class of continuous distributions called the Odds xgamma - G family by taking $W[G(x)]=\frac{G(x;\xi)}{1-G(x;\xi)}$ and $r(t)=\frac{\lambda^2}{1+\lambda}\left(1+\frac{\lambda}{2}t^2\right)e^{-\lambda t}, t > 0, \lambda > 0,$ where $G(x;\xi)$ is a baseline cdf, which depends on a parameter vector $\xi$ and $\bar{G}(x;\xi)= 1- G(x;\xi)$ is the baseline survival function.\\
The distribution function of Odds xgamma - G family of distribution is given by
\begin{eqnarray}\label{eq1}
	F(x;\lambda,\xi)&=& \int_{0}^{\frac{G(x;\xi)}{1-G(x;\xi)}}\frac{\lambda^2}{1+\lambda}\left(1+\frac{\lambda}{2}t^2\right)e^{-\lambda t}dt
	\nonumber\\&=&1-\frac{1+\lambda+\lambda \frac{G(x;\xi)}{\bar{G}(x;\xi)}+\frac{\lambda^2}{2}{\left\{\frac{G(x;\xi)}{\bar{G}(x;\xi)}\right\}}^2}{1+\lambda}e^{-\lambda \frac{G(x;\xi)}{\bar{G}(x;\xi)}}
\end{eqnarray}
The probability density function (pdf) of Odds xgamma - G family of distribution, is given by
\begin{eqnarray}\label{eq2}
	f(x;\lambda,\xi)= \frac{\lambda^2}{1+\lambda}\frac{g(x;\xi)}{[\bar{G}(x;\xi)]^2}\left[1+\frac{\lambda}{2}{\left\{\frac{G(x;\xi)}{\bar{G}(x;\xi)}\right\}}^2\right]e^{-\lambda \frac{G(x;\xi)}{\bar{G}(x;\xi)}}
\end{eqnarray}

The survival function of Odds xgamma - G family of distribution is given by
\begin{eqnarray}\label{eq3}
	S(x;\alpha,\xi)&=&\frac{1+\lambda+\lambda \frac{G(x;\xi)}{\bar{G}(x;\xi)}+\frac{\lambda^2}{2}{\left\{\frac{G(x;\xi)}{\bar{G}(x;\xi)}\right\}}^2}{1+\lambda}e^{-\lambda \frac{G(x;\xi)}{\bar{G}(x;\xi)}}
\end{eqnarray}

The hazard rate function of Odds xgamma - G family of distribution is given by
\begin{eqnarray}\label{eq4}
	h(t;\lambda,\xi)&=&\frac{f(t;\lambda,\xi)}{S(t;\lambda,\xi)}\nonumber\\&=&\frac{\lambda^2 \frac{g(t;\xi)}{[\bar{G}(t;\xi)]^2}\left[1+\frac{\lambda}{2}{\left\{\frac{G(t;\xi)}{\bar{G}(t;\xi)}\right\}}^2\right]}{1+\lambda+\lambda \frac{G(t;\xi)}{\bar{G}(t;\xi)}+\frac{\lambda^2}{2}{\left\{\frac{G(t;\xi)}{\bar{G}(t;\xi)}\right\}}^2}
\end{eqnarray}

\begin{center}
Table 1: Distributions and corresponding $G(x;\xi)/\bar{G}(x;\xi)$ functions
\end{center}
\begin{tabular}{@{} l @{}}
\hline
Distribution~~~~~~~~~~~~~~~~~~~~~~~~~~~~~~~~~~~~~~~~~~~~~~$G(x;\xi)/\bar{G}(x;\xi)$~~~~~~~~~~~~~~~~~~~~~~~~~~~~~~~~~~~~~~~~~~~$\xi$~~~~~~~~\\
\hline
Uniform($0<x<\theta$)~~~~~~~~~~~~~~~~~~~~~~~~~~~~~~~~~~~~~~~$x/(\theta-x)$~~~~~~~~~~~~~~~~~~~~~~~~~~~~~~~~~~~~~~~~~~~~~~~~$\theta$\\
Exponential($x>0$)~~~~~~~~~~~~~~~~~~~~~~~~~~~~~~~~~~~~~~~~$e^{\lambda x}-1$~~~~~~~~~~~~~~~~~~~~~~~~~~~~~~~~~~~~~~~~~~~~~~~~~~$\lambda$\\
Weibull($x>0$)~~~~~~~~~~~~~~~~~~~~~~~~~~~~~~~~~~~~~~~~~~~~~~$e^{\lambda x^{\gamma}}-1$~~~~~~~~~~~~~~~~~~~~~~~~~~~~~~~~~~~~~~~~~~~~~~$(\lambda, \gamma)$\\
Frechet($x>0$)~~~~~~~~~~~~~~~~~~~~~~~~~~~~~~~~~~~~~~~~~~~~~~$(e^{\lambda x^{\gamma}}-1)^{-1}$~~~~~~~~~~~~~~~~~~~~~~~~~~~~~~~~~~~~~~~~~$(\lambda, \gamma)$\\
Half-logistic($x>0$)~~~~~~~~~~~~~~~~~~~~~~~~~~~~~~~~~~~~~~~~$(e^{x}-1)/2$~~~~~~~~~~~~~~~~~~~~~~~~~~~~~~~~~~~~~~~~~~~~~~$\phi$\\
Power function($0<x<1/\theta$)~~~~~~~~~~~~~~~~~~~~~~~~~~~~$[(\theta x)^{-k}-1]^{-1}$~~~~~~~~~~~~~~~~~~~~~~~~~~~~~~~~~~~~~~$(\theta,k)$\\
Pareto($x\geq \theta$)~~~~~~~~~~~~~~~~~~~~~~~~~~~~~~~~~~~~~~~~~~~~~~~$(x/\theta)^k-1$~~~~~~~~~~~~~~~~~~~~~~~~~~~~~~~~~~~~~~~~~~~~$(\theta,k)$\\
Burr XII($>0$)~~~~~~~~~~~~~~~~~~~~~~~~~~~~~~~~~~~~~~~~~~~~~~~$[1+(x/s)^c]^{k}-1 $~~~~~~~~~~~~~~~~~~~~~~~~~~~~~~~~~~$(s,k,c)$\\
Log-logistic($x>0$)~~~~~~~~~~~~~~~~~~~~~~~~~~~~~~~~~~~~~~~~~$[1+(x/s)^c]-1 $~~~~~~~~~~~~~~~~~~~~~~~~~~~~~~~~~~~~~$(s,c)$\\
Lomax($x>0$)~~~~~~~~~~~~~~~~~~~~~~~~~~~~~~~~~~~~~~~~~~~~~~~$[1+(x/s)]^{k}-1 $~~~~~~~~~~~~~~~~~~~~~~~~~~~~~~~~~~~~~$(s,k)$\\
Gumbel($-\infty<x<\infty$)~~~~~~~~~~~~~~~~~~~~~~~~~~~~~~~~$[exp[exp(-(x-\mu)/\sigma)]-1]^{-1}$~~~~~~~~~~~~~~~~~~~~~~$(\mu,\sigma)$\\
Kumaraswamy($0<x<1$)~~~~~~~~~~~~~~~~~~~~~~~~~~~~~~~~$(1-x^a)^{-b}-1$~~~~~~~~~~~~~~~~~~~~~~~~~~~~~~~~~~~~~~$(a,b)$\\
Normal($-\infty<x<\infty$)~~~~~~~~~~~~~~~~~~~~~~~~~~~~~~$\Phi((x-\mu)/\sigma)/(1-\Phi((x-\mu)/\sigma))$~~~~~~~~~~~~~~~~~~~~$(\mu,\sigma)$\\

\hline
\end{tabular}

\section{Some Special Models for Odds xgamma - G Family}
In this section, some new special distributions, namely, Odds xgamma-Uniform, Odds xgamma-Exponential, and GAW-log logistic are introduced.

\subsection{Odds xgamma - Uniform Distribution}
Considering the baseline distribution is Uniform on the interval $(0,\theta),~\theta>0~$with the pdf and cdf, respectively
\begin{eqnarray*}
g(x;\theta)=\frac{1}{\theta}~;0<x<\theta<\infty, ~G(x,\theta)=\frac{x}{\theta}
\end{eqnarray*}
The cdf of Odds xgamma-Uniform distribution is obtained by substituting the pdf and cdf of uniform in $(\ref{eq2})$ as follows
\begin{eqnarray*}
F(x;\lambda,\theta)&=&1-\frac{1+\lambda+\frac{\lambda x}{\theta-x}+\frac{\lambda^2 x^2}{2 (\theta-x)^2}}{1+\lambda}e^{-\frac{\lambda x}{\theta-x}}
\end{eqnarray*}
The corresponding pdf is given by
\begin{eqnarray*}
f(x;\lambda,\theta)&=&\frac{\lambda^2}{1+\lambda}\frac{\theta}{(\theta-x)^2}\left[1+\frac{\lambda x^2}{2(\theta-x)^2}\right]e^{-\frac{\lambda x}{\theta-x}} ~~;0<x<\theta<\infty, \lambda>0.
\end{eqnarray*}
The survival and hazard rate functions are given respectively as follows
\begin{eqnarray*}
R(x;\lambda,\theta)&=&\frac{1+\lambda+\frac{\lambda x}{\theta-x}+\frac{\lambda^2 x^2}{2 (\theta-x)^2}}{1+\lambda}e^{-\frac{\lambda x}{\theta-x}}
\end{eqnarray*}
\begin{eqnarray*}
h(x;\lambda,\theta)&=&\frac{\lambda^2 \theta\left[1+\frac{\lambda x^2}{2(\theta-x)^2}\right]}{(\theta-x)^2 \left[1+\lambda+\frac{\lambda x}{\theta-x}+\frac{\lambda^2 x^2}{2 (\theta-x)^2}\right]}
\end{eqnarray*}

\subsection{Odds xgamma - Exponential Distribution}
Considering the baseline distribution is Exponential with parameter $\theta>0~$. The pdf and cdf are
\begin{eqnarray*}
g(x;\theta)=\theta e^{-\theta x}~;0<x,\theta<\infty, ~G(x,\theta)=1-e^{-\theta x}
\end{eqnarray*}
The cdf of Odds xgamma-Exponential distribution is obtained by substituting the pdf and cdf of uniform in $(\ref{eq2})$ as follows
\begin{eqnarray*}
F(x;\lambda,\theta)&=&1-\frac{1+\lambda e^{\theta x}+\frac{\lambda^2 (e^{\theta x}-1)^2}{2}}{1+\lambda}e^{-\lambda (e^{\theta x}-1)}
\end{eqnarray*}
The corresponding pdf is given by
\begin{eqnarray*}
f(x;\lambda,\theta)&=&\frac{\lambda^2}{1+\lambda}\theta e^{\theta x}\left[1+\frac{\lambda (e^{\theta x}-1)^2}{2}\right]e^{-\lambda (e^{\theta x}-1)} ~~;0<x,\theta<\infty, \lambda>0.
\end{eqnarray*}
The survival and hazard rate functions are given respectively as follows
\begin{eqnarray*}
R(x;\lambda,\theta)&=&\frac{1+\lambda e^{\theta x}+\frac{\lambda^2 (e^{\theta x}-1)^2}{2}}{1+\lambda}e^{-\lambda (e^{\theta x}-1)}
\end{eqnarray*}
\begin{eqnarray*}
h(x;\lambda,\theta)&=&\frac{\lambda^2 \theta e^{\theta x}\left[1+\frac{\lambda (e^{\theta x}-1)^2}{2}\right]}{1+\lambda e^{\theta x}+\frac{\lambda^2 (e^{\theta x}-1)^2}{2}}
\end{eqnarray*}

\subsection{Odds xgamma - Burr XII Distribution}
Considering the baseline distribution is Burr XII (see Burr (1942)) with the following pdf and cdf
\begin{eqnarray*}
g(x;\alpha,\theta)=\alpha\theta x^{(\theta-1)}\left(1+x^\alpha\right)^{-(\theta+1)}~~;x\geq0,\alpha,\theta>0, 
\end{eqnarray*}
\begin{eqnarray*}
G(x;\alpha,\theta)=1- \left(1+x^\alpha\right)^{-\theta}~~;x\geq0,\alpha,\theta>0.
\end{eqnarray*}
The cdf of Odds xgamma-Burr XII distribution is obtained by substituting the pdf and cdf of uniform in $(\ref{eq2})$ as follows
\begin{eqnarray*}
F(x;\lambda,\alpha,\theta)&=&1-\frac{1+\lambda (1+x^{\alpha})^{\theta}+\frac{\lambda^2}{2} \left[(1+x^{\alpha})^{\theta}-1\right]^2}{1+\lambda}e^{-\lambda \left[(1+x^{\alpha})^{\theta}-1\right]}
\end{eqnarray*}
The corresponding pdf is given by
\begin{eqnarray*}
f(x;\lambda,\alpha,\theta)&=&\frac{\lambda^2}{1+\lambda}\alpha \theta x^{\theta-1}(1+x^{\alpha})^{\theta-1}\left[1+\frac{\lambda}{2}\left[(1+x^{\alpha})^{\theta}-1\right]^2 \right]e^{-\lambda \left[(1+x^{\alpha})^{\theta}-1\right]} ~~;0<x,\theta,\alpha<\infty, \lambda>0.
\end{eqnarray*}
The survival and hazard rate functions are given respectively as follows
\begin{eqnarray*}
R(x;\lambda,\alpha,\theta)&=&\frac{1+\lambda (1+x^{\alpha})^{\theta}+\frac{\lambda^2}{2} \left[(1+x^{\alpha})^{\theta}-1\right]^2}{1+\lambda}e^{-\lambda \left[(1+x^{\alpha})^{\theta}-1\right]}
\end{eqnarray*}
\begin{eqnarray*}
h(x;\lambda,\alpha,\theta)&=&\frac{\lambda^2 \alpha \theta x^{\theta-1}(1+x^{\alpha})^{\theta-1}\left[1+\frac{\lambda}{2}\left[(1+x^{\alpha})^{\theta}-1\right]^2 \right]}{1+\lambda (1+x^{\alpha})^{\theta}+\frac{\lambda^2}{2} \left[(1+x^{\alpha})^{\theta}-1\right]^2}
\end{eqnarray*}

\subsection{Odds xgamma - Normal Distribution}
The cdf of Odds xgamma-Normal distribution is obtained by
\begin{eqnarray*}
F(x;\lambda,\mu,\sigma)&=&1-\frac{1+\lambda + \lambda \frac{\Phi(\frac{x-\mu}{\sigma})}{1-\Phi(\frac{x-\mu}{\sigma})}+\frac{\lambda^2}{2}\left[\frac{\Phi(\frac{x-\mu}{\sigma})}{1-\Phi(\frac{x-\mu}{\sigma})}\right]^2}{1+\lambda}e^{-\lambda \frac{\Phi(\frac{x-\mu}{\sigma})}{1-\Phi(\frac{x-\mu}{\sigma})}}
\end{eqnarray*}
The corresponding pdf is given by
\begin{eqnarray*}
f(x;\lambda,\mu,\sigma)&=&\frac{\lambda^2}{(1+\lambda)\sigma}\frac{\phi(\frac{x-\mu}{\sigma})}{[1-\Phi(\frac{x-\mu}{\sigma})]^2}\left[1+\frac{\lambda}{2}\frac{\Phi(\frac{x-\mu}{\sigma})}{1-\Phi(\frac{x-\mu}{\sigma})}\right]^{2} e^{-\lambda \frac{\Phi(\frac{x-\mu}{\sigma})}{1-\Phi(\frac{x-\mu}{\sigma})}} ~~;-\infty<x<\infty.
\end{eqnarray*}

The survival and hazard rate functions are given respectively as follows
\begin{eqnarray*}
R(x;\lambda,\mu,\sigma)&=&\frac{1+\lambda + \lambda \frac{\Phi(\frac{x-\mu}{\sigma})}{1-\Phi(\frac{x-\mu}{\sigma})}+\frac{\lambda^2}{2}\left[\frac{\Phi(\frac{x-\mu}{\sigma})}{1-\Phi(\frac{x-\mu}{\sigma})}\right]^2}{1+\lambda}e^{-\lambda \frac{\Phi(\frac{x-\mu}{\sigma})}{1-\Phi(\frac{x-\mu}{\sigma})}}
\end{eqnarray*}
\begin{eqnarray*}
h(x;\lambda,\mu,\sigma)&=&\frac{\frac{\lambda^2}{\sigma}\frac{\phi(\frac{x-\mu}{\sigma})}{[1-\Phi(\frac{x-\mu}{\sigma})]^2}\left[1+\frac{\lambda}{2}\frac{\Phi(\frac{x-\mu}{\sigma})}{1-\Phi(\frac{x-\mu}{\sigma})}\right]^{2}}{1+\lambda + \lambda \frac{\Phi(\frac{x-\mu}{\sigma})}{1-\Phi(\frac{x-\mu}{\sigma})}+\frac{\lambda^2}{2}\left[\frac{\Phi(\frac{x-\mu}{\sigma})}{1-\Phi(\frac{x-\mu}{\sigma})}\right]^2}
\end{eqnarray*}

\begin{figure}[ht]
\centering
\includegraphics[scale=0.7]{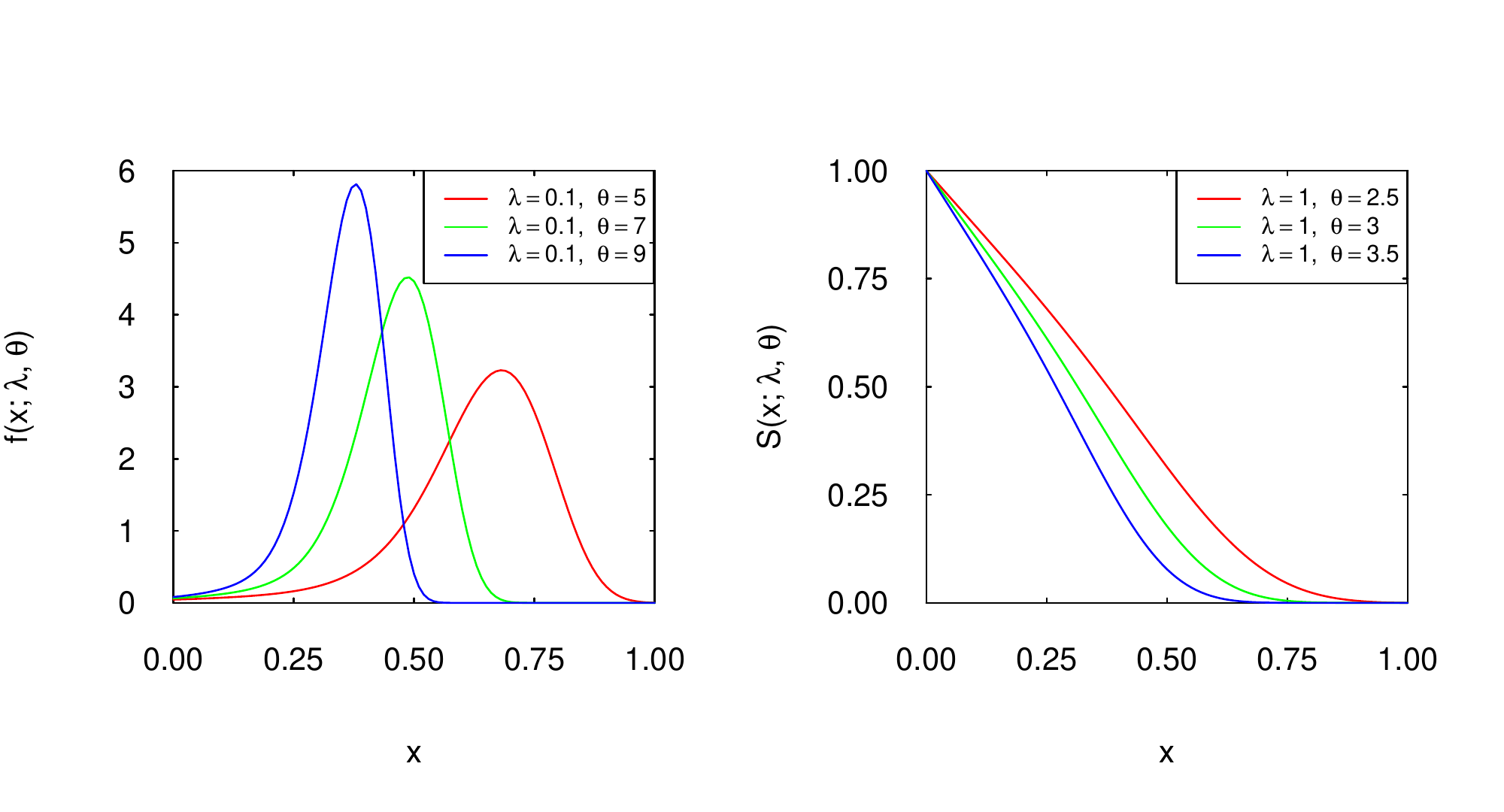} 
\vspace{.1cm} \caption{The pdf and survival function of Odds xgamma - Exponential Distribution}
\label{fig1}
\end{figure}

\begin{figure}[ht]
\centering
\includegraphics[scale=0.7]{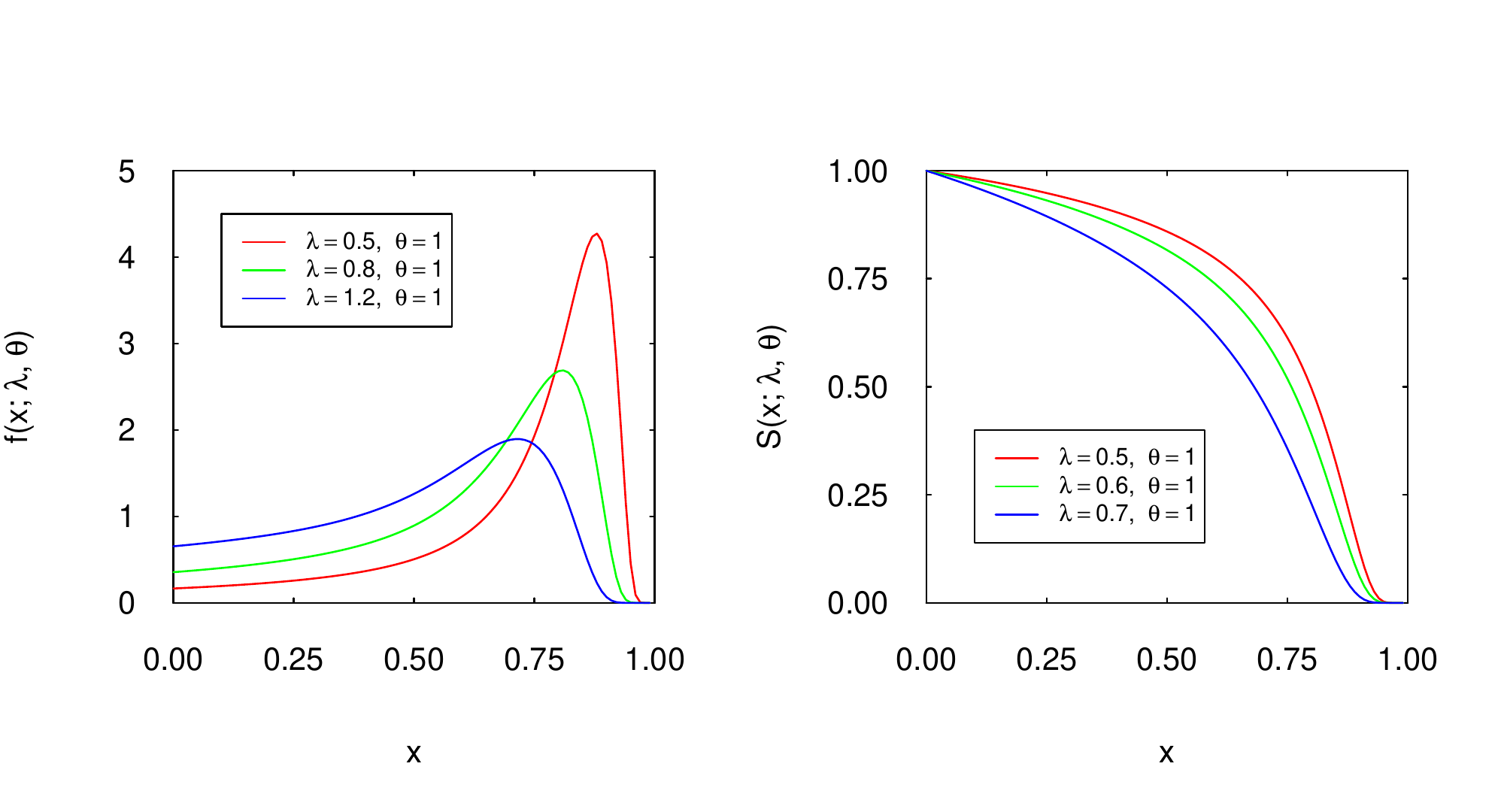} 
\vspace{.1cm} \caption{The pdf and survival function of Odds xgamma - Uniform Distribution}
\label{fig2}
\end{figure}

\section{Some Mathematical Properties}

In this section, some general results of the Odds xgamma - G family are derived.
\subsection{Mixture Representation}
Expansion formulae of the Odds xgamma - G family, such as; the pdf and cdf are derived.
The probability density function (pdf) of Odds xgamma - G family of distribution, is given by
\begin{eqnarray}\label{eq10}
f(x;\lambda,\xi)&=&\frac{\lambda^2}{1+\lambda}\frac{g(x;\xi)}{[\bar{G}(x;\xi)]^2}\left[1+\frac{\lambda}{2}{\left\{\frac{G(x;\xi)}{\bar{G}(x;\xi)}\right\}}^2\right]e^{-\lambda\frac{G(x;\xi)}{\bar{G}(x;\xi)}}\nonumber\\&=&\frac{\lambda^2}{1+\lambda}\frac{g(x;\xi)}{[\bar{G}(x;\xi)]^2}\left[1+\frac{\lambda}{2}{\left\{\frac{G(x;\xi)}{\bar{G}(x;\xi)}\right\}}^2\right]\sum_{i=0}^\infty\frac{(-1)^i}{i!}{\lambda^i}\left[\frac{G(x;\xi)}{\bar{G}(x;\xi)}\right]^i\nonumber\\&=&\sum_{i=0}^\infty\frac{(-1)^i}{i!}\frac{\lambda^{i+2}}{1+\lambda}\frac{g(x;\xi)[G(x;\xi)]^i}{[\bar{G}(x;\xi)]^{i+2}}\left[1+\frac{\lambda}{2}{\left\{\frac{G(x;\xi)}{\bar{G}(x;\xi)}\right\}}^2\right]\nonumber\\&=&\sum_{i=0}^\infty\frac{(-1)^i}{i!}\frac{\lambda^{i+2}}{1+\lambda}g(x;\xi)[G(x;\xi)]^i[\bar{G}(x;\xi)]^{-(i+2)}+\sum_{i=0}^\infty\frac{(-1)^i}{i!}\frac{\lambda^{i+3}}{2(1+\lambda)}g(x;\xi)[G(x;\xi)]^{i+2}[\bar{G}(x;\xi)]^{-(i+4)}\nonumber\\&=&\sum_{i,j=0}^\infty\frac{(-1)^i}{i!}\binom{i+j+1}{j}\frac{\lambda^{i+2}}{1+\lambda}g(x;\xi)[G(x;\xi)]^{i+j}+\sum_{i,k=0}^\infty\frac{(-1)^i}{i!}\binom{i+k+3}{k}\frac{\lambda^{i+3}}{2(1+\lambda)}g(x;\xi)[G(x;\xi)]^{i+k+2}\nonumber\\&=&\sum_{i,j=0}^\infty w_{ij} h_{i+j+1}(x;\xi)+\sum_{i,k=0}^\infty w_{ik} h_{i+k+3}(x;\xi)
\end{eqnarray}
Where, $w_{ij}=\frac{(-1)^i}{i!}\binom{i+j+1}{j}\frac{\lambda^{i+2}}{1+\lambda}$, $w_{ik}=\frac{(-1)^i}{i!}\binom{i+k+3}{k}\frac{\lambda^{i+3}}{2(1+\lambda)}$ and $ h_{m+n+1}(x;\xi)=g(x;\xi)[G(x;\xi)]^{m+n}$

The cdf of X can be given by integrating equation $(\ref{eq10})$ as
\begin{eqnarray}\label{eq11}
F(x;\lambda,\xi)&=&\sum_{i,j=0}^\infty w_{ij} H_{i+j+1}(x;\xi)+\sum_{ik=0}^\infty w_{i,k} H_{i+k+3}(x;\xi)
\end{eqnarray}

\subsection{Shapes of the Odds xgamma - G family of distribution}
The shapes of the density and hazard rate functions can also be described analytically.\\
Now, 
\begin{eqnarray*}
	f(x;\lambda,\xi)= \frac{\lambda^2}{1+\lambda}\frac{g(x;\xi)}{[\bar{G}(x;\xi)]^2}\left[1+\frac{\lambda}{2}{\left\{\frac{G(x;\xi)}{\bar{G}(x;\xi)}\right\}}^2\right]e^{-\lambda \frac{G(x;\xi)}{\bar{G}(x;\xi)}}
\end{eqnarray*}
So, \begin{eqnarray*}	\ln f(x;\lambda,\xi)= \ln \left(\frac{\lambda^2}{1+\lambda}\right)+ \ln g(x;\xi) - 2\ln \bar{G}(x;\xi) + \ln \left[1+\frac{\lambda}{2}{\left\{\frac{G(x;\xi)}{\bar{G}(x;\xi)}\right\}}^2\right] -\lambda \frac{G(x;\xi)}{\bar{G}(x;\xi)}
\end{eqnarray*}
Now, The critical points of the Odds xgamma - G density function are the roots of the equation:
\begin{eqnarray*}	\frac{d}{dx}\ln f(x;\lambda,\xi)= \frac{g'(x;\xi)}{g(x;\xi)}+ \frac{2 g(x;\xi)}{\bar{G}(x;\xi)}+\frac{\frac{\lambda g(x;\xi) G(x;\xi) }{[\bar{G}(x;\xi)]^3}}{1+\frac{\lambda}{2}{\left\{\frac{G(x;\xi)}{\bar{G}(x;\xi)}\right\}}^2}-\frac{\lambda g(x;\xi)}{[\bar{G}(x;\xi)]^2}=0
\end{eqnarray*}

\subsection{Quantile function}
The quantile function, say $Q(u)=F^{-1}(u)$,  of the Odds xgamma - G family is derived by inverting $(\ref{eq1})$ as follows
\begin{eqnarray*}
	u&=&1-\frac{1+\lambda+ \frac{\lambda Q(u)}{1-Q(u)}+\frac{\lambda^2}{2}{\left\{\frac{Q(u)}{1-Q(u)}\right\}}^2}{1+\lambda}e^{-\lambda \frac{Q(u)}{1-Q(u)}}
\end{eqnarray*}
So,
\begin{eqnarray*}
	1+\lambda+ \frac{\lambda Q(u)}{1-Q(u)}+\frac{\lambda^2}{2}{\left\{\frac{Q(u)}{1-Q(u)}\right\}}^2&=&(1+\lambda)(1-u)e^{\lambda \frac{Q(u)}{1-Q(u)}}
\end{eqnarray*}
Taking Log on both sides, after some simplifications, the previous equation is reduced to
\begin{eqnarray}\label{eq7}
	 \frac{\lambda Q(u)}{1-Q(u)}+\ln (1+\lambda)+ \ln (1-u) -\ln\left[1+\lambda+\lambda \frac{Q(u)}{1-Q(u)}+\frac{\lambda^2}{2}{\left\{\frac{Q(u)}{1-Q(u)}\right\}}^2\right]&=&0
\end{eqnarray}
By solving the nonlinear equation $(\ref{eq7})$, numerically, the Odds xgamma - G family random variable X can be generated, where u has the uniform distribution on the unit interval.

\subsection{Moments}
The rth moment of random variable X can be obtained from pdf $(\ref{eq10})$ as follows
\begin{eqnarray*}
	\mu_r^{'}&=&\int_0^{\infty}x^rf(x,\xi)dx\nonumber\\&=& \sum_{i,j=0}^\infty w_{ij}\int_0^{\infty}x^r h_{i+j+1}(x;\xi)dx+ \sum_{i,k=0}^\infty w_{ik}\int_0^{\infty}x^r h_{i+k+3}(x;\xi)dx
\end{eqnarray*}
Therefore
\begin{eqnarray}\label{eq21}
	\mu_r^{'}&=& \sum_{i,j=0}^{\infty}w_{ij}I_{i,j,r}+\sum_{i,k}^{\infty}w_{ik}I_{i,k,r};~r=1,~2,~....
\end{eqnarray}
where,
$I_{i,j,r}=\int_0^{\infty}x^r h_{i+j+1}(x;\xi)dx$ and $I_{i,k,r}= \int_0^{\infty}x^r h_{i+k+3}(x;\xi)dx.$\\
In particular, the mean and variance of Odds xgamma - G family are obtained as follows:
\begin{eqnarray*}
E(X)&=& \sum_{i,j=0}^{\infty}w_{ij}I_{i,j,1}+\sum_{i,k=0}^{\infty}w_{ik}I_{i,k,1}
\end{eqnarray*}
The variance is
\begin{eqnarray*}
Var(X)&=&\sum_{i,j=0}^{\infty}w_{ij}I_{i,j,2}+\sum_{i,k=0}^{\infty}w_{ik}I_{i,k,2}-\left[\sum_{i,j=0}^{\infty}w_{ij}I_{i,j,1}+\sum_{i,k=0}^{\infty}w_{ik}I_{i,k,1}\right]^2
\end{eqnarray*}
Additionally, measures of skewness and kurtosis of family can be obtained, based on $(\ref{eq21})$, according to the following relations
\begin{eqnarray*}
	\gamma_{1}&=&\frac{\mu_3^{'}-3\mu_2^{'}\mu_1^{'}+2\mu_1^{'^3}}{\left(\mu_2^{'}-\mu_1^{'^2}\right)^{3/2}}
\end{eqnarray*}
\begin{eqnarray*}
	\gamma_{2}&=&\frac{\mu_4^{'}-4\mu_3^{'}\mu_1^{'}+6\mu_2^{'}\mu_1^{'^2}-3\mu_1^{'^4}}{\left(\mu_2^{'}-\mu_1^{'^2}\right)^2}
\end{eqnarray*}

\subsection{Generating Function}
\textbf{The Moment Generating function(MGF)} of Odds xgamma - G family is defined as
\begin{eqnarray*}
M_{X}(t)&=&\sum_{r=0}^\infty \frac{t^{r}}{r!} \mu_{r}^{'}
\end{eqnarray*}
where, $\mu_r^{'}$  is the $r^{th}$ moment about origin, then  the moment generating function of Odds xgamma - G family is obtained by using $(\ref{eq21})$ as follows
\begin{eqnarray*}
M_{X}(t)&=&\sum_{r=0}^\infty \frac{t^{r}}{r!} \left[\sum_{i,j=0}^{\infty}w_{ij}I_{i,j,r}+\sum_{i,k}^{\infty}w_{ik}I_{i,k,r}\right]
\end{eqnarray*}

\textbf{Characteristic Function(CF)}: \begin{eqnarray*} \Psi_{X}(t)&=& E(e^{itX})\nonumber\\&=&\sum_{r=0}^\infty \frac{(it)^{r}}{r!} \mu_{r}^{'}\nonumber\\&=&\sum_{r=0}^\infty\frac{(it)^{r}}{r!}\left[\sum_{i,j=0}^{\infty}w_{ij}I_{i,j,r}+\sum_{i,k}^{\infty}w_{ik}I_{i,k,r}\right]\end{eqnarray*}

\textbf{Cumulant Generating Function(CGF)}: \begin{eqnarray*} K_{X}(t)&=&\ln_{e}(M_{X}(t))\nonumber\\&=&\ln_{e}\sum_{r=0}^\infty \frac{t^{r}}{r!} \left[\sum_{i,j=0}^{\infty}w_{ij}I_{i,j,r}+\sum_{i,k}^{\infty}w_{ik}I_{i,k,r}\right]\end{eqnarray*}

\subsection{Entropy}
An entropy of a random variable X is a measure of variation of the uncertainty. A popular entropy measure is Renyi entropy (Renyi 1961). If X has the probability density function f(x), then Renyi entropy is defined by \begin{eqnarray}\label{eq15}H_R(\beta)&=&\frac{1}{1-\beta}\ln\left\{\int_{a}^{\infty}f^\beta(x)dx \right\}\end{eqnarray}
where $\beta>0$ and $\beta \neq 1$. \\
Here, we derive expressions for the Renyi entropy for the Odds xgamma-G distribution. Due to the fact that the parameter is not in general a natural number, it is difficult to use $(\ref{eq21})$ for entropy derivation. So, we use $(\ref{eq21})$, the power series for the exponential and the generalized binomial expansion to obtain the Renyi entropy of X.

Now \begin{eqnarray*} \int_{0}^{\infty}f^\beta(x)dx &=&\left[\frac{\lambda^2}{1+\lambda}\right]^{\beta}\int_{0}^{\infty}\left[\frac{g(x;\xi)}{[\bar{G}(x;\xi)]^2}\right]^{\beta}\left[1+\frac{\lambda}{2}{\left\{\frac{G(x;\xi)}{\bar{G}(x;\xi)}\right\}}^2\right]^{\beta}e^{-\lambda \beta \frac{G(x;\xi)}{\bar{G}(x;\xi)}}dx\nonumber\\&=&\left[\frac{\lambda^2}{1+\lambda}\right]^{\beta}\int_{0}^{\infty}\left[\frac{g(x;\xi)}{[\bar{G}(x;\xi)]^2}\right]^{\beta}\left[1+\frac{\lambda}{2}{\left\{\frac{G(x;\xi)}{\bar{G}(x;\xi)}\right\}}^2\right]^{\beta}\sum_{i=0}^\infty\frac{(-1)^i}{i!}{\lambda ^i}{\beta ^i}\left[\frac{G(x;\xi)}{\bar{G}(x;\xi)}\right]^{i}dx\nonumber\\&=&\left[\frac{\lambda^2}{1+\lambda}\right]^{\beta}\sum_{i=0}^\infty\frac{(-1)^i}{i!}{\lambda ^i}{\beta^i}\int_{0}^{\infty}\left[\frac{g(x;\xi)}{[\bar{G}(x;\xi)]^2}\right]^{\beta}\left[1+\frac{\lambda}{2}{\left\{\frac{G(x;\xi)}{\bar{G}(x;\xi)}\right\}}^2\right]^{\beta}\left[\frac{G(x;\xi)}{\bar{G}(x;\xi)}\right]^{i}dx\nonumber\\&=&\left[\frac{\lambda^2}{1+\lambda}\right]^{\beta}\sum_{i=0}^\infty\frac{(-1)^i}{i!}{\lambda^i}{\beta^i}\int_{0}^{\infty}\left[\frac{g(x;\xi)}{[\bar{G}(x;\xi)]^2}\right]^{\beta}\left[\frac{G(x;\xi)}{\bar{G}(x;\xi)}\right]^{i}\sum_{j=0}^\beta \binom{\beta}{j}\left[\frac{\lambda}{2}\right]^{j}\left[\frac{G(x;\xi)}{\bar{G}(x;\xi)}\right]^{2j}dx\nonumber\\&=&\left[\frac{\lambda^2}{1+\lambda}\right]^{\beta}\sum_{i=0}^\infty\sum_{j=0}^\beta\frac{(-1)^i}{i!}{\lambda^i}{\beta^i}\binom{\beta}{j}\left[\frac{\lambda}{2}\right]^{j}\int_{0}^{\infty}\left[g(x;\xi)\right]^{\beta}\left[G(x;\xi)\right]^{i+2j}\left[\bar{G}(x;\xi)\right]^{-(i+2j+2\beta)}dx\nonumber\\&=&\left[\frac{\lambda^2}{1+\lambda}\right]^{\beta}\sum_{i=0}^\infty\sum_{j=0}^\beta\sum_{k=0}^\infty\frac{(-1)^i{\lambda^{i+j}}{\beta^i}}{i!2^j}\binom{\beta}{j}\binom{i+2j+2\beta-1}{k}\int_{0}^{\infty}\left[g(x;\xi)\right]^{\beta}\left[G(x;\xi)\right]^{i+2j+k}dx\nonumber\\&=&\left[\frac{\lambda^2}{1+\lambda}\right]^{\beta}\sum_{i=0}^\infty\sum_{j=0}^\beta\sum_{k=0}^\infty\frac{(-1)^i{\lambda^{i+j}}{\beta^i}}{i!2^j}\binom{\beta}{j}\binom{i+2j+2\beta-1}{k}K(\beta,i,j,k)\end{eqnarray*}

\subsection{Order Statistics}
A branch of statistics known as order statistics plays a proeminent role in real-life applications involving data relating to life testing studies. These statistics are required in many fields, such as climatology, engineering and industry, among others. A comprehensive exposition of order statistics and associated inference is provided by David and Nagaraja (2003). Let $X_{r:n}$ denote the $r^{th}$ order statistic. The density $f_{r:n}(x)$ of the $r^{th}$ order statistic, for $r = 1(1)n,$ from independent and identically distributed random variables $X_1, X_2,.....  X_n$ having the Odds xgamma-G distribution is given by
\begin{eqnarray*}
f_{r:n}(x)&=&M\left[F(x)\right]^{r-1}\left[1-F(x)\right]^{n-r}f(x).
\end{eqnarray*}
\begin{eqnarray*}
f_{r:n}(x)&=&M\sum_{s=0}^{n-r}(-1)^s\binom{n-r}{s}\left[F(x)\right]^{r+s-1}f(x).
\end{eqnarray*}
Where $M=\frac{n!}{(r-1)!(n-r)!}$

\begin{eqnarray}\label{eq22}
f_{r:n}(x;\Phi)&=& M\sum_{s=0}^{n-r}(-1)^s\binom{n-r}{s}\left[1-\frac{1+\lambda+\lambda \frac{G(x;\xi)}{\bar{G}(x;\xi)}+\frac{\lambda^2}{2}{\left\{\frac{G(x;\xi)}{\bar{G}(x;\xi)}\right\}}^2}{1+\lambda}e^{-\lambda \frac{G(x;\xi)}{\bar{G}(x;\xi)}}\right]^{r+s-1}\nonumber\\&&.\frac{\lambda^2}{1+\lambda}\frac{g(x;\xi)}{[\bar{G}(x;\xi)]^2}\left[1+\frac{\lambda}{2}{\left\{\frac{G(x;\xi)}{\bar{G}(x;\xi)}\right\}}^2\right]e^{-\lambda \frac{G(x;\xi)}{\bar{G}(x;\xi)}}\nonumber\\&=& M \frac{\lambda^2}{1+\lambda}\sum_{s=0}^{n-r}(-1)^s\binom{n-r}{s}\frac{g(x;\xi)}{[\bar{G}(x;\xi)]^2}\left[1+\frac{\lambda}{2}{\left\{\frac{G(x;\xi)}{\bar{G}(x;\xi)}\right\}}^2\right]e^{-\lambda \frac{G(x;\xi)}{\bar{G}(x;\xi)}}\nonumber\\&&.\sum_{k=0}^{r+s-1}(-1)^k\binom{r+s-1}{k}\left[\frac{1+\lambda+\lambda \frac{G(x;\xi)}{\bar{G}(x;\xi)}+\frac{\lambda^2}{2}{\left\{\frac{G(x;\xi)}{\bar{G}(x;\xi)}\right\}}^2}{1+\lambda}\right]^{k}e^{-k\lambda \frac{G(x;\xi)}{\bar{G}(x;\xi)}}\nonumber\\&=& M \frac{\lambda^2}{1+\lambda}\sum_{s=0}^{n-r}\sum_{k=0}^{r+s-1}\sum_{i=0}^{\infty}(-1)^{s+k+i}\frac{[\lambda(k+1)]^{i}}{i!}\binom{n-r}{s}\binom{r+s-1}{k}\frac{g(x;\xi)[G(x;\xi)]^i}{[\bar{G}(x;\xi)]^{i+2}}\nonumber\\&&.\left[\frac{1+\lambda+\lambda \frac{G(x;\xi)}{\bar{G}(x;\xi)}+\frac{\lambda^2}{2}{\left\{\frac{G(x;\xi)}{\bar{G}(x;\xi)}\right\}}^2}{1+\lambda}\right]^{k}\left[1+\frac{\lambda}{2}{\left\{\frac{G(x;\xi)}{\bar{G}(x;\xi)}\right\}}^2\right]
\end{eqnarray}

\subsection{Stress-Strength Reliability}

The measure of reliability of industrial components has many applications especially in the area of engineering. The reliability of a product (system) is the probability that the product (system) will perform its intended function for a specified time period when operating under normal (or stated) environmental conditions. The component fails at the instant that the random stress $X_2$ applied to it exceeds the random strength $X_1$, and the component will function satisfactorily whenever $X_1 > X_2$. Hence, $R = P(X_2 < X_1)$ is a measure of component reliability (see Kotz, Lai, and Xie (2003)). We derive the reliability R when $X_1$ and $X_2$ have independent Oxgamma-G(x; $\lambda_1;\xi$) and Oxgamma-G(x; $\lambda_2;\xi$) distributions with the same parameter vector $\xi$ for the baseline G. The reliability is denoted by
\begin{eqnarray*} R \nonumber&=&\int_{0}^{\infty}f_{1}(x)F_{2}(x)dx\end{eqnarray*}
The pdf of $X_1$ and cdf of $X_2$ are obtained from equation $(\ref{eq10})$ and $(\ref{eq11})$ as 
\begin{eqnarray*}
f_{1}(x)&=&\sum_{i,j=0}^\infty p_{ij}(\lambda_1) g(x,\xi)[G(x,\xi)]^{i+j}+\sum_{i,k=0}^\infty p_{ik}(\lambda_1) g(x,\xi)[G(x,\xi)]^{i+k+2}
\end{eqnarray*}
\begin{eqnarray*}
F_{2}(x)&=&\sum_{l,m=0}^\infty q_{lm}(\lambda_2)[G(x,\xi)]^{l+m}+\sum_{l,n=0}^\infty q_{ln}(\lambda_2) [G(x,\xi)]^{l+n+2}
\end{eqnarray*}
Where \begin{eqnarray*} p_{ij}(\lambda_1)=\frac{(-1)^i}{i!}\binom{i+j+1}{j}\frac{\lambda^{i+2}}{1+\lambda}\end{eqnarray*}
\begin{eqnarray*} p_{ik}(\lambda_1)=\frac{(-1)^i}{i!}\binom{i+k+3}{k}\frac{\lambda^{i+3}}{2(1+\lambda)}\end{eqnarray*}
\begin{eqnarray*} q_{lm}(\lambda_2)=\frac{(-1)^l}{l!}\binom{l+m+1}{m}\frac{\lambda^{l+2}}{1+\lambda}\end{eqnarray*}
\begin{eqnarray*} q_{ln}(\lambda_2)=\frac{(-1)^l}{l!}\binom{l+n+3}{n}\frac{\lambda^{l+3}}{2(1+\lambda)}\end{eqnarray*}

Hence, \begin{eqnarray*} R \nonumber&=&\sum_{i,j,l,m=0}^\infty p_{ij}(\lambda_1)q_{lm}(\lambda_2)\int_{0}^{\infty}g(x,\xi)[G(x,\xi)]^{i+j+l+m}dx+\nonumber\\&& \sum_{i,j,l,n=0}^\infty p_{ij}(\lambda_1)q_{ln}(\lambda_2)\int_{0}^{\infty}g(x,\xi)[G(x,\xi)]^{i+j+l+n+2}dx+\nonumber\\&& \sum_{i,k,l,m=0}^\infty p_{ik}(\lambda_1)q_{lm}(\lambda_2)\int_{0}^{\infty}g(x,\xi)[G(x,\xi)]^{i+k+l+m+2}dx+\nonumber\\&& \sum_{i,k,l,n=0}^\infty p_{ik}(\lambda_1)q_{ln}(\lambda_2)\int_{0}^{\infty}g(x,\xi)[G(x,\xi)]^{i+k+l+n+4}dx\end{eqnarray*}

Setting $u=G(x,\xi)$, the reliability of the Odds xgamma - G family distribution reducess to
\begin{eqnarray*} R \nonumber&=&\sum_{i,j,l,m=0}^\infty\frac{p_{ij}(\lambda_1)q_{lm}(\lambda_2)}{i+j+l+m+1} +\sum_{i,j,l,n=0}^\infty \frac{p_{ij}(\lambda_1)q_{ln}(\lambda_2)}{i+j+l+n+3}\nonumber\\&&+ \sum_{i,k,l,m=0}^\infty\frac{p_{ik}(\lambda_1)q_{lm}(\lambda_2)}{i+k+l+m+3} + \sum_{i,k,l,n=0}^\infty \frac{p_{ik}(\lambda_1)q_{ln}(\lambda_2)}{i+k+l+n+5}\end{eqnarray*}

\begin{figure}[ht]
\centering
\includegraphics[scale=0.7]{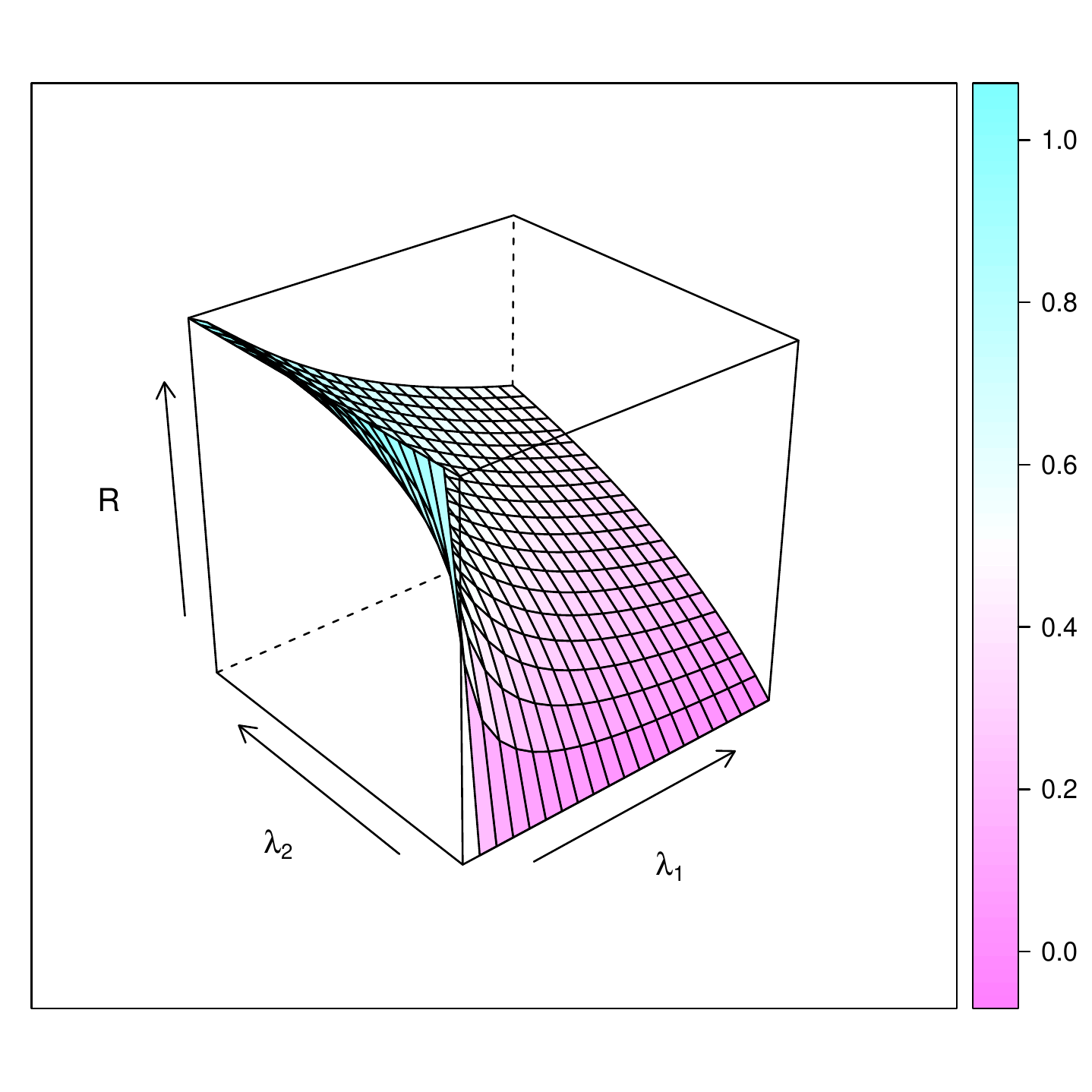} 
\vspace{.1cm} \caption{Stress-Strength Reliability, R for different $\lambda_1$ and $\lambda_2$ when $\theta_1=\theta_2$, of odds xgamma Exponential Distribution}
\label{fig3}
\end{figure}

\subsection{Incomplete Moments, Mean Deviations and Lorenz and Benferroni Curves}

The r-th incomplete moment, say,$m_r^I(t)$, of the GAW-G distribution is given by
\begin{eqnarray*}
	m_r^I(t)&=& \int_0^tx^rf(x,\Phi)dx.
\end{eqnarray*}
We can write from equation ($\ref{eq18}$),
\begin{eqnarray}\label{eq4.4.1}
	m_r^I(t)&=& \int_0^tx^r\left[\sum_{i,j=0}^{\infty}w_{ij}g(x,\xi)[G(x,\xi)]^{i+j}+\sum_{i,k=0}^{\infty}w_{ik}g(x,\xi)[G(x,\xi)]^{i+k+2}\right]dx.
\end{eqnarray}
\begin{ex}
Consider the Odds- xgamma uniform distribution discussed in subsection 2.1. 
\begin{eqnarray*}
	m_r^I(t)&=& \sum_{i,j=0}^{\infty}w_{ij}\frac{t^{i+j+r+1}}{\theta^{i+j+1}(i+j+r+1)}+\sum_{i,k=0}^{\infty}w_{ik}\frac{t^{i+k+r+3}}{\theta^{i+k+3}(i+k+r+3)}
\end{eqnarray*}
\end{ex}

The amount of scatter in a population is evidently measured to some extent by the totality of the deviations from the mean and median. The mean deviations about the mean $\delta_1=E(\mid X-\mu_1^{'}\mid)$ and median $\delta_2=E(\mid X-M\mid)$ of $X$ can be used as measures of spread (or dispersion) in a population besides range and standard deviation. They are given by $\delta_1=2\mu_1^{'}F(\mu_1^{'})-2m_1^I(\mu_1^{'})$ and $\delta_2=\mu_1^{'}-2m_1^I(M)$,respectively. Here, $\mu_1^{'}=E(X)$ is to be obtained from ($\ref{eq21}$) with $r=1$, $F(\mu_1^{'})$ is to calculated from ($\ref{eq1}$), $m_1^I(\mu_1^{'})$ is the first incomplete function obtained from ($\ref{eq4.4.1}$) with $r=1$ and $M$ is the median of $X$ obtained by solving ($\ref{eq7}$) for $u=0.5$.\\

The Lorenz and Benferroni curves are defined by $L(p)=m_1^I(x_p)/\mu_1^{'}$ and $B(p)=m_1^I(x_p)/(p\mu_1^{'})$, respectively, where $x_p=F^{-1}(p)$ can be computed numerically by ($\ref{eq7}$) with $u=p$. These curves have significant roles in economics, reliability, demography, insurance and medicine. For details in this aspect, the readers are referred to Pundir et al.(2005) and references cited therein.
\begin{figure}[ht]
\centering
\includegraphics[scale=0.7]{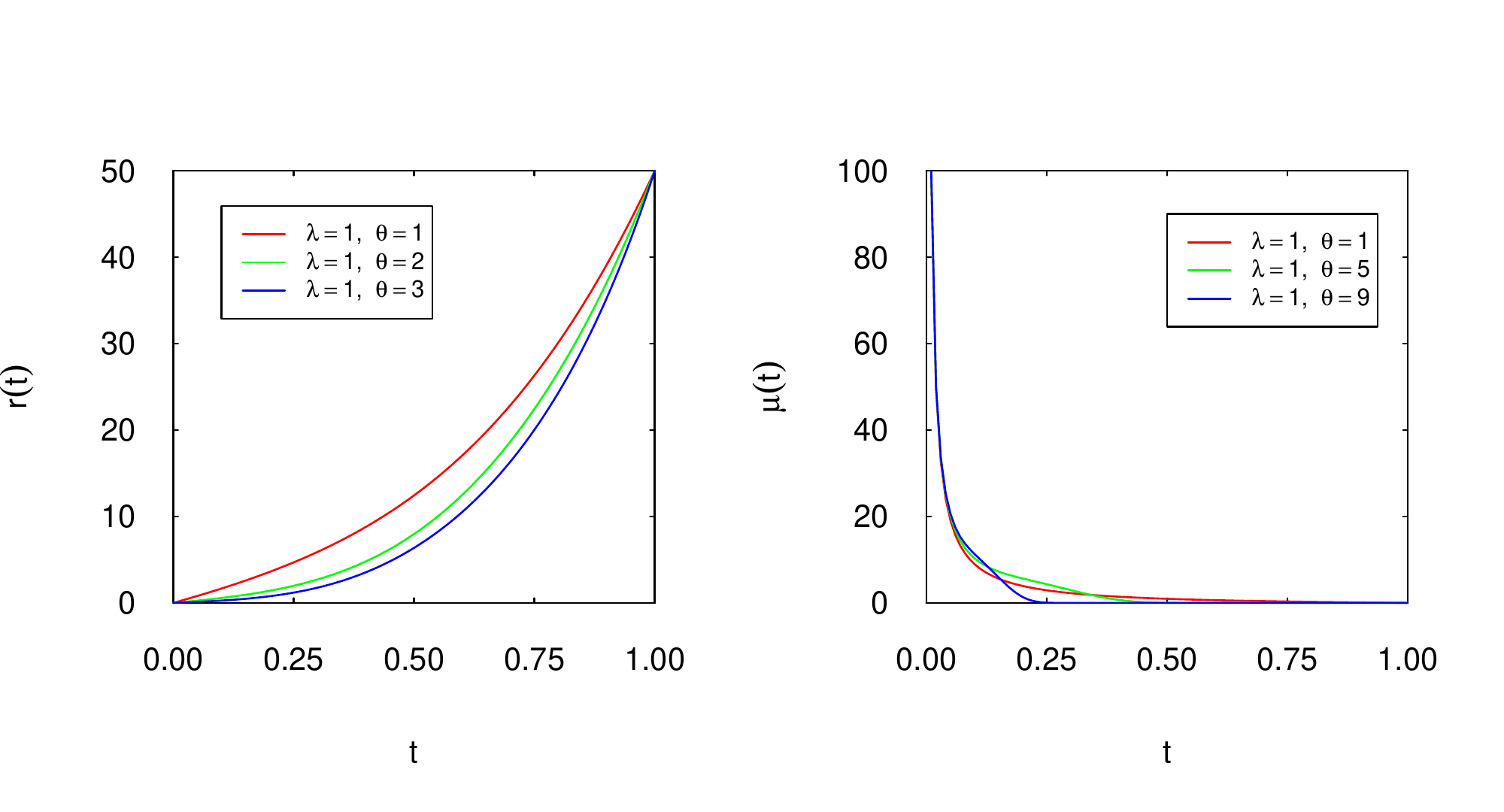} 
\vspace{.1cm} \caption{The Hazard Rate and Reversed Hazard Rate of Odds xgamma - Exponential Distribution}
\label{fig4}
\end{figure}

\begin{figure}[ht]
\centering
\includegraphics[scale=0.7]{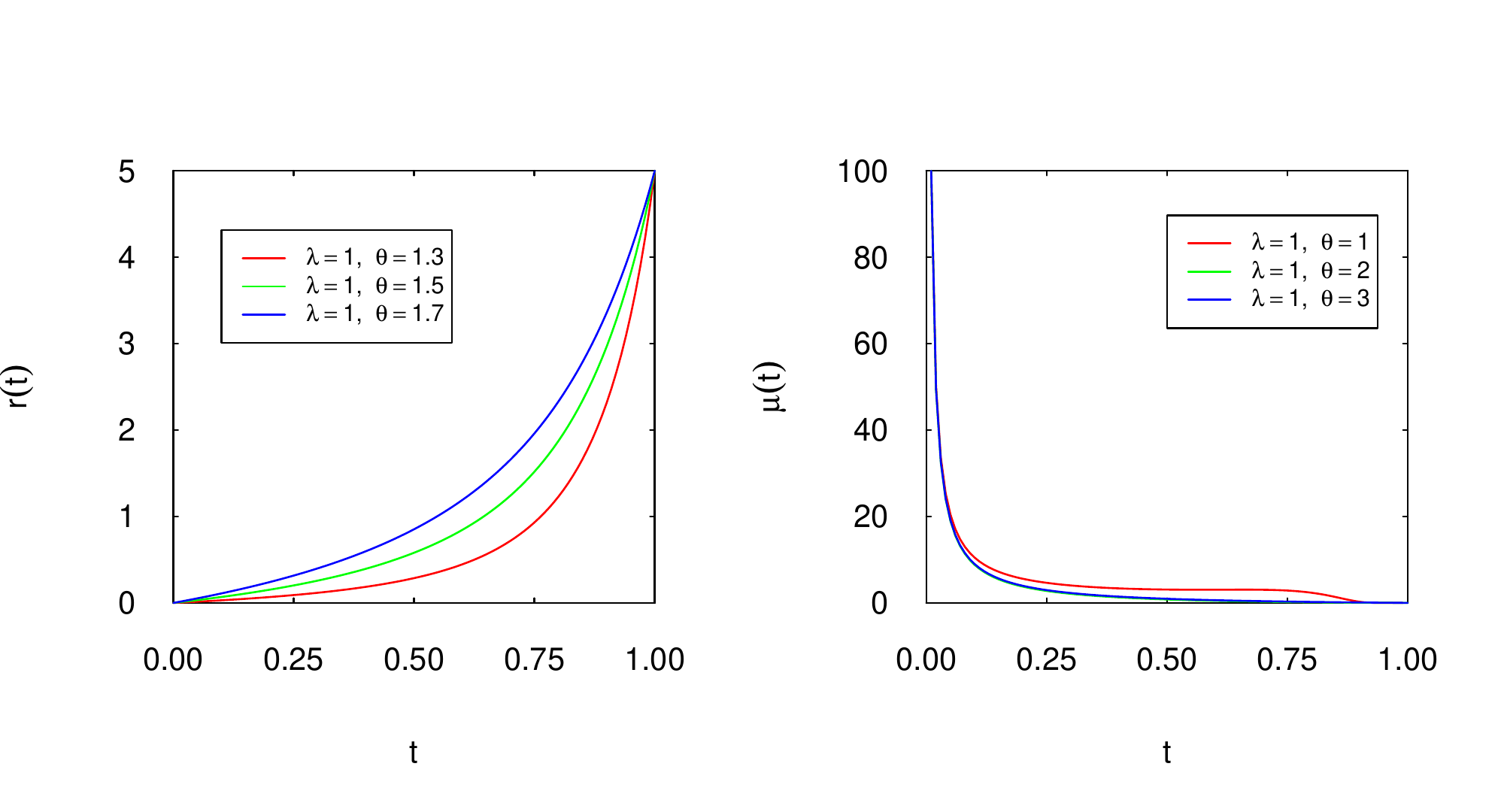} 
\vspace{.1cm} \caption{The Hazard Rate and Reversed Hazard Rate of Odds xgamma - Uniform Distribution}
\label{fig5}
\end{figure}

\subsection{Moments of the residual life}
The hazard rate, mean residual life, left truncated mean function are some functions related to the residual lifetime of an unit. These functions uniquely determine the cumulative distribution function (cdf), $F(x)$. See, for instance, Gupta(1975) and Zoroa et al.(1990).
\begin{d1}
Let $X$, be a random variable denoting the lifetime of a unit is at age $t$. Then $X_t=X-t\mid X>t$ denotes the remaining lifetime beyond that age $t$.
\end{d1}
The cdf $F(x)$ is uniquely determined by the r-th moment of the residual life of $X$ (for $r=1,~2,~...$)[Navarro et al.(1998)], and it is given by
\begin{eqnarray*}
	m_r(t)=E[X_t]&=& \frac{1}{\bar F(t)}\int_t^{\infty}(x-t)^rdF(x)\\&=&\frac{1}{1-F(t)}\int_t^{\infty}(x-t)^rf(x,\Phi)dx
\end{eqnarray*}
In particular, if $r=1$, then $m_1(t)$ represents an interesting function called the mean residual life (MRL) function that indicates the expected life length for a unit which is alive at age $t$. The MRL function has wide spectrum of applications in reliability/survival analysis, social studies, biomedical sciences, economics, population study, insurance industry, maintenance and product quality control and product technology.
\begin{ex}
Consider again the Odds- xgamma uniform distribution discussed in subsection 2.1. 
\begin{eqnarray*}
\bar{F}(t)&=&\frac{1+\lambda+\frac{\lambda t}{\theta-t}+\frac{\lambda^2 t^2}{2 (\theta-t)^2}}{1+\lambda}e^{-\frac{\lambda t}{\theta-t}}
\end{eqnarray*}
Using ($\ref{eq10}$), we have
\begin{eqnarray*}
	\int_t^{\theta}(x-t)^rf(x,\Phi)dx&=& \sum_{i,j=0}^{\infty}\frac{w_{ij}}{\theta^{i+j+1}}\sum_{u=0}^r\binom{r}{u}(-t)^u\frac{\theta^{i+j+r-u+1}-t^{i+j+r-u+1}}{i+j+r-u+1}\\&&+\sum_{i,k=0}^{\infty}\frac{w_{ik}}{\theta^{i+k+3}}\sum_{u=0}^r\binom{r}{u}(-t)^u\frac{\theta^{i+k+r-u+3}-t^{i+k+r-u+3}}{i+k+r-u+3}
\end{eqnarray*}
For the MRL function,
 \begin{eqnarray*}
	\int_t^{\theta}(x-t)f(x,\Phi)dx&=& \sum_{i,j=0}^{\infty}\frac{w_{ij}}{\theta^{i+j+1}}\left[\frac{\theta^{i+j+2}-t^{i+j+2}}{i+j+2}-\frac{t\left(\theta^{i+j+1}-t^{i+j+1}\right)}{i+j+1}\right]\\&&+\sum_{i,k=0}^{\infty}\frac{w_{ik}}{\theta^{i+k+3}}\left[\frac{\theta^{i+k+3}-t^{i+k+3}}{i+k+3}-\frac{t\left(\theta^{i+k+2}-t^{i+k+2}\right)}{i+k+2}\right]
\end{eqnarray*}
\end{ex}
\subsection{Moments of the reversed residual life}
In many real life situations uncertainty is not necessarily related to the future but can also refer to the past. For instance, consider a system whose state is observed only at certain preassigned inspection time $t$. If the system is inspected for the first time and it is found to be "`down"', then failure relies on the past (i.e. on which instant in $(0,t)$ it has failed). Therefore, study of a notion that is dual to the residual life, in the sense that it refers to the past time and not to future seems worthwhile [see Di Crescenzo and Longobardi (2002)].
\begin{d1}
Let $X$ be a random variable denoting the lifetime of a unit is down at age $t$. Then $\bar{X}_t=t-X\mid X<t$ denotes the idle time or inactivity time or reversed residual life of the unit at age $t$.
\end{d1}
In case of forensic science, people may be interested in estimating $\bar{X}_t$ in order to ascertain the exact time of death of a person. In Insurance industry, it represents the period remained unpaid by a policy holder due to death. For details, see Block et al.(1998), Chandra and Roy(2001), Maiti and Nanda(2009), and Nanda et al.(2003).
The r-th moment of $\bar{X}_t$ (for $r=1,~2,~...$) is given by
\begin{eqnarray*}
	\bar{m}_r(t)=E[\bar{X}_t]&=& \frac{1}{F(t)}\int_0^{t}(t-x)^rdF(x)\\&=&\frac{1}{F(t)}\int_0^{t}(t-x)^rf(x,\Phi)dx
\end{eqnarray*}
In particular, if $r=1$, then $\bar{m}_1(t)$ represents a function called the mean idle time or inactivity time (MIT) or reversed residual life (MRRL) function that indicates the expected inactive life length for a unit which is first observed down at age $t$. The properties of MIT function have been explored by Ahmad et al. (2005) and Kayid and Ahmad (2004).
\begin{ex}
Consider again the Odds- xgamma uniform distribution discussed in subsection 2.1.
\begin{eqnarray*}
F(t)&=&1-\frac{1+\lambda+\frac{\lambda t}{\theta-t}+\frac{\lambda^2 t^2}{2 (\theta-t)^2}}{1+\lambda}e^{-\frac{\lambda t}{\theta-t}}
\end{eqnarray*}
Using ($\ref{eq10}$), we have
\begin{eqnarray*}
	\int_0^{t}(t-x)^rf(x,\Phi)dx&=& \sum_{i,j=0}^{\infty}\frac{w_{ij}}{\theta^{i+j+1}}\sum_{u=0}^r(-1)^u\binom{r}{u}\frac{t^{i+j+r+1}}{i+j+u+1}\\&&+\sum_{i,k=0}^{\infty}\frac{w_{ik}}{\theta^{i+k+3}}\sum_{u=0}^r(-1)^u\binom{r}{u}\frac{t^{i+k+r+3}}{i+k+u+3}
\end{eqnarray*}
For the MIT (or MRRL) function,
 \begin{eqnarray*}
	\int_0^{t}(t-x)f(x,\Phi)dx&=& \sum_{i,j=0}^{\infty}\frac{w_{ij}}{\theta^{i+j+1}}.\frac{t^{i+j+2}}{(i+j+1)(i+j+2)}\\&&+\sum_{i,k=0}^{\infty}\frac{w_{ik}}{\theta^{i+k+3}}.\frac{t^{i+k+4}}{(i+k+3)(i+k+4)}
\end{eqnarray*}
\end{ex}

\begin{figure}[ht]
\centering
\includegraphics[scale=0.7]{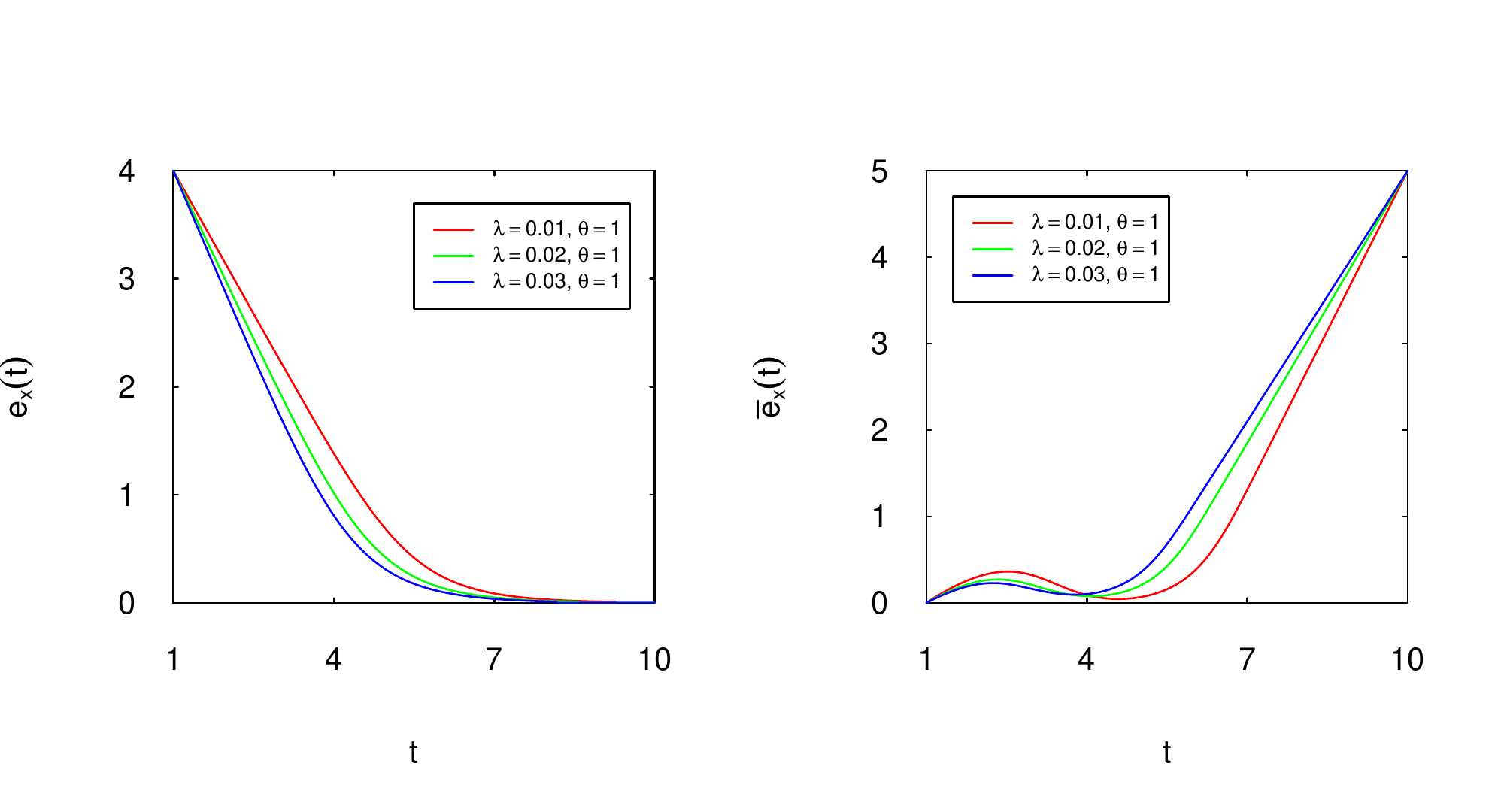} 
\vspace{.1cm} \caption{Mean Residual Life and Reversed Mean Residual Life  of the Odds xgamma - Exponential Distribution}
\label{fig6}
\end{figure}

\section{Maximum Likelihood Estimation}
In this section, we determine the Maximum likelihood estimates(MLEs) of the model parameters of the new family of distributions from complete samples only. Let$ x_1, x_2,....,x_n $be a observed values from the Odds xgamma -G family distribution with parameters $\lambda$ and $\xi$. Let $\Phi=(\lambda, \xi)^{T}$ be the p x 1 parameter vector. The total log-likelihood function for $\Phi$ is given by
\begin{eqnarray*}
l(\Phi)&=&2n\ln\lambda-n\ln(1+\lambda)+ \sum_{i=1}^{n}\ln g(x_i;\xi) -2\sum_{i=1}^{n}\ln\bar{G}(x_i;\xi)+ \sum_{i=1}^{n}\ln\left[1+\frac{\lambda}{2}{\left\{\frac{G(x_i;\xi)}{\bar{G}(x_i;\xi)}\right\}}^2\right]-\lambda \sum_{i=1}^{n}\frac{G(x_i;\xi)}{\bar{G}(x_i;\xi)}\nonumber\\&=&2n\ln\lambda-n\ln(1+\lambda)+ \sum_{i=1}^{n}\ln g(x_i;\xi)  -2\sum_{i=1}^{n}\ln\bar{G}(x_i;\xi)+ \sum_{i=1}^{n}\ln\left[1+\frac{\lambda}{2}\left[V(x_i;\xi)\right]^2\right]-\lambda \sum_{i=1}^{n}V(x_i;\xi)
\end{eqnarray*}

where $V(x_i;\xi)=G(x_i;\xi)/\bar{G}(x_i;\xi)$. The components of the score function $U\left(\Phi\right)=\left(U_{\lambda}, U_{\xi}\right)^{T}$ are

\begin{eqnarray*}
U_{\lambda}=\frac{2n}{\lambda}-\frac{n}{1+\lambda}-\sum_{i=1}^{n}V(x_i;\xi)
\end{eqnarray*}
and
\begin{eqnarray*}
U_{\xi_{k}}=\sum_{i=1}^{n}\frac{\frac{\partial}{\partial \xi_{k}}g(x_i;\xi)}{g(x_i;\xi)}-2\sum_{i=1}^{n}\frac{\frac{\partial}{\partial \xi_{k}} \bar{G}(x_i;\xi)}{\bar{G}(x_i;\xi)}+\sum_{i=1}^{n}\frac{\lambda V(x_i;\xi)\frac{\partial}{\partial \xi_{k}}V(x_i;\xi)}{1+\frac{\lambda}{2}\left[V(x_i;\xi)\right]^2}-\lambda \sum_{i=1}^{n} \frac{\partial}{\partial \xi_{k}}V(x_i;\xi)
\end{eqnarray*}

Setting $U_{\lambda}$ and $U_{\xi}$ equal to zero and solving the equations simultaneously yields the MLE $\hat{\Phi}=\left(\hat{\lambda},\hat{\xi}\right)^{T}$ of $\Phi=\left(\lambda,\xi\right)^{T}$. These equations cannot be solved analytically and statistical software can be used to solve them numerically using iterative methods such as the Newton- Raphson type algorithms.

\section{Application}
~~In this section, we fit the Odds xgamma - Exponential Distribution model to two real data sets. The first data set obtained from Smith and Naylor (1987) [24]. The data are the strengths of 1.5 cm glass fibres, measured at the National Physical Laboratory, England and have been shown in Table 1. Histogram shows that the data set is negatively skewed. Frank Gomes-Silva et al. (2017) fitted this data to the Odd Lindley Weibull Distribution. We have fitted this data set with the Odds xgamma - Exponential Distribution. The estimated values of the parameters were $\lambda= 0.08736933$, $\theta= 2.191986$, log-likelihood =$ -14.04618$ and AIC = $32.09237$. Histogram and fitted Odds xgamma - Exponential curve to data have been shown in Figure $\ref{fig14}$.\\

\begin{center}
Table 1: Strengths of glass fibres data set
\end{center}
\begin{tabular}{@{} l @{}}
\hline
0.55 0.93 1.25 1.36 1.49 1.52 1.58 1.61 1.64 1.68 1.73 1.81 2.00 0.74 1.04 1.27 1.39 1.49 1.53 1.59 1.61\\
1.66 1.68 1.76 1.82 2.01 0.77 1.11 1.28 1.42 1.50 1.54 1.60 1.62 1.66 1.69 1.76 1.84 2.24 0.81 1.13 1.29\\
1.48 1.50 1.55 1.61 1.62 1.66 1.70 1.77 1.84 0.84 1.24 1.30 1.48 1.51 1.55 1.61 1.63 1.67 1.70 1.78 1.89\\
\hline
\end{tabular}
\\
\\
\begin{center}
Table 2: Summarized results of fitting different distributions to data set of Smith and Naylor (1987)
\end{center}
\begin{tabular}{@{} l @{}}
\hline
Distribution~~~~~~~~~~~~~~~~~~~~~~~~Estimate of the parameter~~~~~~~~~~~~~~~~~~Log-likelihood~~~~~~~~~~~~AIC\\
\hline
Odd Lindley Weibull~~~~~~~~~~~~~$\hat{a}=0.049, \hat{\alpha}=1.102, \hat{\lambda}=0.492$~~~~~~~~~~~~~~~$-14.193$~~~~~~~~~~~~~~~$34.387$\\
O-xg Exp~~~~~~~~~~~~~~~~~~~~~~~~~~~$\hat{\lambda}=0.08736933, \hat{\theta}=2.191986$~~~~~~~~~~~~~~~~~~$-14.046$~~~~~~~~~~~~~~~$32.092$\\
\hline
\end{tabular}
\\
\\

~~The second data set taken from the R base package. It is located in the Indometh object. The data consists of plasma concentrations of indomethacin(mcg/ml) and have been shown in Table 1. Histogram shows that the data set is positively skewed. Cordeiro and Barreto-Souza (2009) fitted this data to the Beta gamma Distribution. We have fitted this data set with the Odds xgamma - Exponential Distribution. The estimated values of the parameters were $\hat{\lambda}=16.80083, \hat{\theta}=0.1050095$ and AIC = $66.68347$. Histogram and fitted Odds xgamma - Exponential curve to data have been shown in Figure $\ref{fig8}$.\\

\begin{center}
Table 3: Plasma concentrations of indomethacin data set
\end{center}
\begin{tabular}{@{} l @{}}
\hline
1.50~ 0.94~ 0.78~ 0.48~ 0.37~ 0.19~ 0.12~ 0.11~ 0.08~ 0.07~ 0.05~ 2.03~ 1.63~ 0.71~ 0.70~ 0.64~ 0.36~ 0.32\\
0.20~ 0.25~ 0.12~ 0.08~ 2.72~ 1.49~ 1.16~ 0.80~ 0.80~ 0.39~ 0.22~ 0.12~ 0.11~ 0.08~ 0.08~ 1.85~ 1.39~ 1.02\\
0.89~ 0.59~ 0.40~ 0.16~ 0.11~ 0.10~ 0.07~ 0.07~ 2.05~ 1.04~ 0.81~ 0.39~ 0.30~ 0.23~ 0.13~ 0.11~ 0.08~ 0.10\\
0.06~ 2.31~ 1.44~ 1.03~ 0.84~ 0.64~ 0.42~ 0.24~ 0.17~ 0.13~ 0.10~ 0.09\\\hline
\end{tabular}
\\

\begin{center}
Table 4: Summarized results of fitting different distributions for Plasma concentrations of indomethacin data set
\end{center}
\begin{tabular}{@{} l @{}}
\hline
Distribution~~~~~~~~~~~~~~~~~~~~~~~~Estimate of the parameter~~~~~~~~~~~~~~~~~Log-likelihood~~~~~~~~~~~~~~~~~~~AIC\\
\hline
Beta gamma~~~~~~~~~~~~~~~$\hat{a}=0.977, \hat{b}=9.181, \hat{\rho}=1.055, \hat{\lambda}=5.553$~~~~~~~~~~~~-31.368~~~~~~~~~~~~~~~~~~~$70.736$\\
O-xg Exp~~~~~~~~~~~~~~~~~~~~~~~~~$\hat{\lambda}=16.80083, \hat{\theta}=0.1050095$~~~~~~~~~~~~~~~~~~~~~~~-31.341~~~~~~~~~~~~~~~~~~~$66.683$\\
\hline
\end{tabular}

\begin{figure}[ht]
\centering
\includegraphics[scale=0.7]{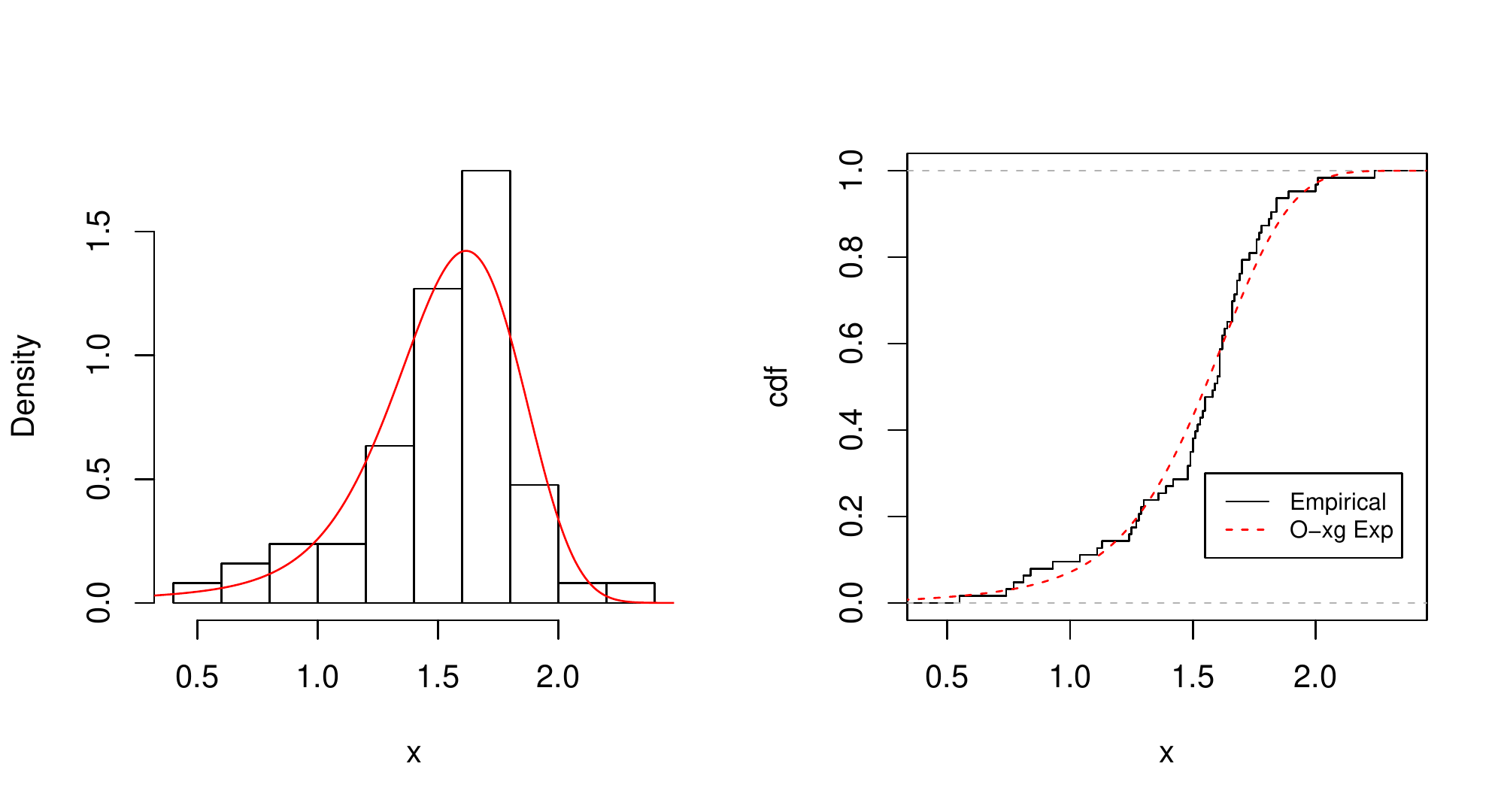} 
\vspace{.1cm} \caption{Plots of the fitted pdf and the estimated cdf of the Odds xgamma - Exponential model for 1.5 cm glass fibres data set}
\label{fig7}
\end{figure}

\begin{figure}[ht]
\centering
\includegraphics[scale=0.7]{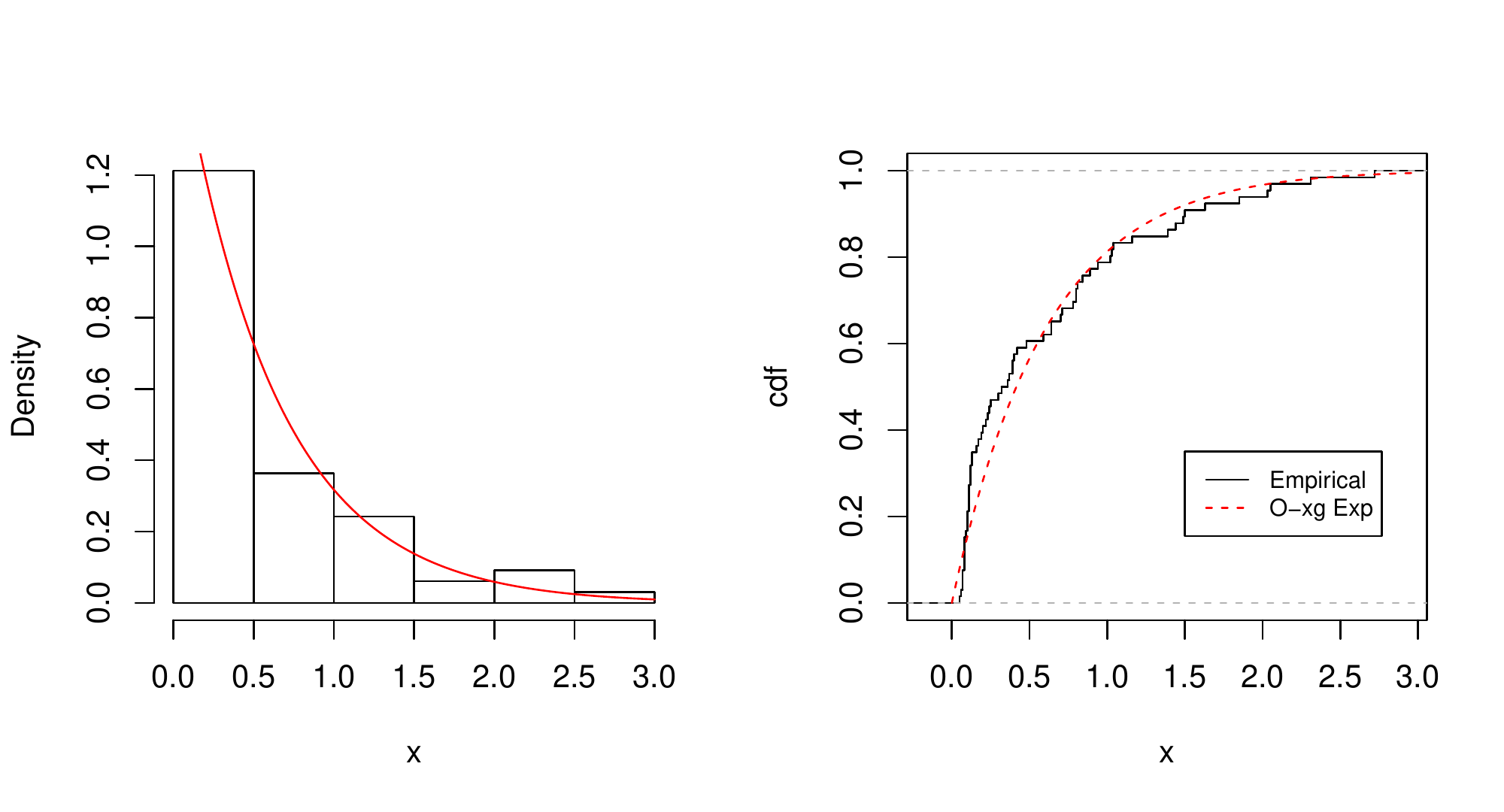} 
\vspace{.1cm} \caption{Plots of the fitted pdf and the estimated cdf of the Odds xgamma - Exponential model for Plasma concentrations of indomethicin data set}
\label{fig8}
\end{figure}

\section{Concluding Remarks}
We have introduced and studied a new generalized family of distributions, called the Odds xgamma - G Family of Distribution. Properties of the Odds xgamma - G Family of Distribution include: an expansion for the density function and expressions for the quantile function, moment generating function, ordinary moments, incomplete moments, mean deviations, Lorenz and Benferroni curves, reliability properties including mean residual life and mean inactivity time, and order statistics. The maximum likelihood method is employed to estimate the model parameters.. Two real data set is used to demonstrate the flexibility of distribution belonging to the introduced family. These special models give better fits than other models. We hope the findings of the paper will be quite useful for the practitioners in various fields of probability, statistics and applied sciences.

\section*{Bibliography}
\begin{enumerate}
\item Alexander C, Cordeiro GM, Ortega EMM, Sarabia JM (2012). "Generalized Beta-generated Distributions." Computational Statistics and Data Analysis, 56, 1880-1897.
\item Alizadeh M, Emadi M, Doostparast M, Cordeiro GM, Ortega EMM, Pescim RR (2015)."A New Family of Distributions: The Kumaraswamy Odd Log-logistic, Properties and Applications." Hacet. J. Math Stat, \textbf{forthcoming}.
\item Alzaatreh A, Lee C, Famoye F (2013). "A New Method for Generating Families of Distributions." Metron, 71, 63-79.
\item Alzaghal A, Felix F, Carl L (2013). "Exponentiated T-X Family of Distributions with Some Applications." International Journal of Statistics and Probability, 2, 31-49.
\item Asgharzedah A, Bakouch HS, Esmaeli H (2013). "Pareto Poisson-Lindley Distribution with Applications." Journal of Applied Statistics, 40, 1717-1734.
\item Bourguignon M, Silva RB, Cordeiro GM (2014). "The Weibull-G Family of Probability Distributions." Journal of Data Science, 12, 53-68.
\item Chandler KN (1952). "The Distribution and Frequency of Record Values." J. Roy. Statist. Soc., Ser B, 14, 220-228.
\item Cooray K (2006). "Generalization of the Weibull Distribution: The Odd Weibull Family." Statistical Modelling, 6, 265-277.
\item Cordeiro GM, Alizadeh M, Ortega EMM (2014). "The Exponentiated Half-logistic Family of Distributions: Properties and Applications." Journal of Probability and Statistics, 2014, 21 pages.
\item Cordeiro GM, de Castro M (2011). "A New Family of Generalized Distributions." Journal of Statistical Computation and Simulation, 81, 883-893.
\item Cordeiro GM, Nadarajah S (2011). "Closed-form Expressions for Moments of a Class of Beta Generalized Distributions." Brazilian Journal Probability and Statistics, 52, 14-33.
\item David HA, Nagaraja HN (2003). Order Statistics. John Wiley \& Sons, New Jersey.
\item Dziubdziela W, Kopocinski B (1976). "Limiting Properties of the $k^{th}$ Sth Record Value." Appl. Math., 15, 187-190.
\item Eugene N, Lee C, Famoye F (2002). "Beta-normal Distribution and Its Applications." Com-munications in Statistics Theory and Methods, 31, 497-512.
\item Exton H (1978). Handbook of Hypergeometric Integrals: Theory, Applications, Tables, Com-puter Programs. New York: Ellis Horwood, New York.
\item Ghitany ME, Atieh B, Nadarajah S (2008). "Lindley Distribution and Its Applications." Mathematics and Computers in Simulation, 78, 493-506.
\item Gradshteyn IS, Ryzhik IM (2007). Table of Integrals, Series, and Products. Academic Press, New York.
\item Gupta P, Singh B (2012). "Parameter Estimation of Lindley Distribution with Hybrid Censored Data." International Journal of System Assurance Engineering and Management, 4, 378-3751.
\item Gupta RC, Gupta RD, Gupta PL (1998). "Modeling Failure Time Data by Lehman Alternatives." Communications in Statistics Theory and Methods, 27, 887-904.
\item Gupta RD, Kundu D (1999). "Generalized Exponential Distributions." Australian and New Zealand Journal of Statistics, 41, 173-188.
\item Johnson NL, Kotz S, Balakrishnan N (1994). Continuous Univariate Distributions. John Wiley, New York.
\item Johnson NL, Kotz S, Balakrishnan N (1995). Continuous Univariate Distributions. John Wiley and Sons, New York.
\item Kotz S, Lai CD, Xie M (2003). "On the Effect of Redundancy for Systems with Dependent Components." IIE Trans, 35, 1103-1110.
\item Lindley DV (1958). "Fiducial Distributions and Bayes' Theorem." Journal of the Royal Statistical Society. Series B. (Methodological), 20, 102-107.
\item Marshall AN, Olkin I (1997). "A New Method for Adding a Parameter to a Family of Distributions with Applications to the Exponential and Weibull Families." Biometrika, 84, 641-552.
\item Mazucheli J, Achcar JA (2011). "The Lindley Distribution Applied to Competing Risks Lifetime Data." Computer Methods and Programs in Biomedicine, 104, 188-192.
\item Mudholkar GS, Srivastava DK (1993). "Exponentiated Weibull Family for Analyzing Bathtub Failure-rate Data." IEEE Transactions on Reliability, 42, 299-302.
\item Mudholkar GS, Srivastava DK, Freimer M (1995). "The Exponentiated Weibull Family: A Reanalysis of the Bus-motor-failure Data." Technometrics, 37, 436-445.
\item Nadarajah Sand Bakouch HS, Tahmasbi R (2011). "A Generalized Lindley Distribution." Sankhya B, 73, 331-359.
\item Nadarajah S (2005). "The Exponentiated Gumbel Distribution with Climate Application." Environmetrics, 17, 13-23.
\item Nadarajah S, Gupta AK (2007). "The Exponentiated Gamma Distribution with Application to Drought Data." Calcutta Statistical Association Bulletin, 59, 29-54.
\item Nadarajah S, Kotz S (2006). "The Exponentiated-type Distributions." Acta Applicandae Mathematicae, 92, 97-111.
\item Renyi A (1961). On Measures of Entropy and Information. In: Proceedings of the 4th Berkeley Symposium on Mathematical, Statistics and Probability. University of California Press, Berkeley.
\item Sankaran M (1970). "The Discrete Poisson-Lindley Distribution." Biometrics, 26, 145-149.
\item Shannon CE (1951). "Prediction and Entropy of Printed English." The Bell System Technical Journal, 30, 50-64.
\item Shirke DT, Kakade CS (2006). "On Exponentiated Lognormal Distribution." International Journal of Agricultural and Statistical Sciences, 2, 319-326.
\item Tahir MH, Cordeiro GM, Alzaatreh A, Mansoor M, Zubair M (2016a). "The Logistic-X Family of Distributions and Its Applications." Commun. Stat. Theory Methods, (forthcoming).
\item Tahir MH, Zubair M, Mansoor M, Cordeiro GM, Alizadeh M, Hamedani GG (2016b). "A New Weibull-G Family of Distributions." Hacet. J. Math. Stat., (forthcoming).
\item Trott M (2006). The Mathematica Guidebook for Symbolics. Springer, New York.
\item Warahena-Liyanage G, Pararai M (2014). "A Generalized Power Lindley Distribution with Applications." Asian Journal of Mathematics and Applications, pp. 1-23.
\item Zakerzadeh Y, Dolati A (2009). "Generalized Lindley Distribution." Journal of Mathematical Extension, 3, 13-25.
\item Zimmer WJ, Keats JB, Wang FK (1998). "The Burr XII Distribution in Reliability Analysis." J. Qual. Technol, 30, 386-394.
\item Zografos K, Balakrishnan N (2009). "On Families of Beta-and Generalized Gamma-generated Distribution and Associate Inference." Statistical Methodology, 6, 344-362.
\end{enumerate}
\end{document}